Article

# Sufficient Conditions for Some Stochastic Orders of Discrete Random Variables with Applications in Reliability

Félix Belzunce [1,*], Carolina Martínez-Riquelme [1] and Magdalena Pereda [2]

[1] Departamento de Estadística e Investigación Operativa, Facultad de Matemáticas, Universidad de Murcia, 30100 Espinardo, Murcia, Spain; carolina.martinez7@um.es
[2] Collège Sciences et Technologies pour l'Energie et l'Environnement (STEE), Université de Pau et des Pays de L'Adour, Avenue de l'Université, CEDEX, BP576-64012 Pau, France; magdalena.peredav@um.es
* Correspondence: belzunce@um.es

**Abstract:** In this paper we focus on providing sufficient conditions for some well-known stochastic orders in reliability but dealing with the discrete versions of them, filling a gap in the literature. In particular, we find conditions based on the unimodality of the likelihood ratio for the comparison in some stochastic orders of two discrete random variables. These results have interest in comparing discrete random variables because the sufficient conditions are easy to check when there are no closed expressions for the survival functions, which occurs in many cases. In addition, the results are applied to compare several parametric families of discrete distributions.

## 1. Introduction

The comparison of random quantities in terms of the so called "stochastic orders" has received great attention along the last 50 years and we can find applications of this topic in different fields such as reliability, insurance, risks, finance, epidemics, and so on (see [1–6]). This topic deals with several criteria to compare two random quantities in order to select which one is "larger", according to their magnitude, dispersion, residual lifetimes, concentration, and so on. An interesting topic in this context is that of finding sufficient conditions for several stochastic orders when there are no closed expressions for some of the functions used in the comparison. Results in this direction can be traced from [7] up to most recent ones by [8–10] and the references therein. Whereas most of the results on this topic are given for continuous random variables, it is hard to find results in the discrete case. As far as we know, there is just one paper related to the comparison of two discrete random variables in the usual stochastic order (see [11]). The purpose of this paper is to fill the gap and provide tools that are easy to deal with, to check if two discrete random variables are ordered according to the most prominent stochastic orders in reliability and survival analysis, which are the hazard rate and mean residual life orders.

Usually, the lifetime of a unit or a living organism is treated as a continuous random variable. However, there are situations where the random lifetime has to be treated as a discrete random variable. For example, in some cases, the number of operations completed correctly prior to the failure of the system or the number of times switching a device to turn it on till it fails is more interesting than the time until failure Moreover, even if the failure time has interest, we deal with discrete times whenever the unit or the system is checked only once per time period. Along the last two decades, there is an increasing interest in the study of reliability modeling and survival analysis in discrete time as can be seen in [12]. Therefore, it is natural to study the hazard rate and mean residual life orders for discrete random variables.

Despite being a powerful tool to compare the residual lifetime of two random lifetimes, the checking of these criteria requires the computation of series which usually do not have a closed expression, which is an important disadvantage in real problems. It is well-known that a sufficient condition for both criteria is the likelihood ratio order, which is satisfied when the ratio of the mass probability (density) functions is monotone. However, there are many situations where this ratio is not monotone but it is unimodal. In the continuous case, the unimodality of the likelihood ratio has been proved to be a sufficient condition, under some additional mild conditions, for the hazard rate and mean residual life orders. In the discrete case, ref. [11] related the unimodality of the ratio of two mass probability functions to the usual stochastic order, as well as applying the results to compare some well-known discrete families of distributions.

In this paper, we show that there are situations where the unimodality of the likelihood ratio is a sufficient condition for the hazard rate order (which is stronger than the stochastic order) or the mean residual life order. These results will be applied to the comparison of some of the most important discrete families of distributions in the context of reliability in discrete time.

The organization of the paper is as follows. In Section 2 we state the notation and main definitions of reliability measures and stochastic orders that are considered along the paper. The main results are provided in Section 3, where the unimodality of the likelihood ratio is related to the comparison of hazard rate and mean residual life functions. Later, in Section 4, we provide applications of these results to the comparison of specific discrete parametric models of interest in reliability.

## 2. Notations and Definitions

In this section we provide the definition of the main reliability measures and stochastic orders considered in the context of reliability in discrete time. In order to provide a unified framework, we have followed [12]. Along this paper "increasing" and "decreasing" mean "nondecreasing" and "nonincreasing," respectively. For any random variable $X$ and any event $A$, we use $\{X|A\}$ to denote the random variable whose distribution is the conditional distribution of $X$ given $A$. Given a discrete random variable $X$ the mass probability function will be denoted by $f_X(x) = P(X = x)$, the distribution function by $F_X(x) = P(X \leq x)$ and the support by $Supp(X)$. The lower and upper extremes of the support will be denoted by $l_X$ and $u_X$, respectively. Given that in this context a random variable $X$ represents the random lifetime of an item or a living organism, the random variable is non negative with $Supp(X) \subseteq \mathbb{N}_0 = \mathbb{N} \cup \{0\}$ and $l_X = 0$.

First, we provide the definitions of the survival, hazard rate, and mean residual life functions.

**Definition 1.** *Let X be a discrete random variable with mass probability function $f_X$.*
*We define the following reliability measures:*

(i) *The survival function of X, denoted by $S_X$, is defined as*

$$S_X(x) = P(X \geq x) = \sum_{i=x}^{u_X} f_X(i), \text{ for all } x \in \mathbb{N}_0.$$

*This function provides the probability of being alive at time x for a unit with random lifetime X.*

(ii) *The hazard rate function of X, denoted by $h_X$, is defined as*

$$h_X(x) = P(X = x | X \geq x) = \frac{f_X(x)}{S_X(x)}, \text{ for all } x \in \mathbb{N}_0 \text{ such that } S_X(x) > 0.$$

*The hazard rate at time X is the probability to fail at time x for a unit that has survived up to time x, which is not the case for continuous distributions. Therefore, it is a local measure.*

(iii) The mean residual life function of X, denoted by $m_X$, is defined as

$$m_X(x) = E[X - x | X > x] = \frac{\sum_{i=x+1}^{\infty} S_X(i)}{S_X(x+1)}, \text{ for all } x \in \mathbb{N}_0 \cup \{-1\}$$

such that $S_X(x+1) > 0$.

The mean residual life is the expected additional time for a unit that has survived beyond x. This function is more informative than the hazard rate, given that it takes into account the history of the unit beyond x.

**Remark 1.** *(i) If the support of X is $\{0, 1, \ldots, u_X\}$, with $u_X < \infty$, then $h_X(u_X) = 1$. As a convention, we consider $h_X(x) = 1$, for all $x > u_X$. It is clear that $0 \leq h_X(x) \leq 1$, for all $x \in \mathbb{N}_0$. Therefore, the interpretation and the boundedness of the hazard rate is different from that of the continuous case.*

*(ii) From Definition 1 (iii) trivially follows that $m_X(-1) = E[X] + 1$.*

Now, we introduce the discrete versions of the stochastic orders based on the previous reliability measures. Additionally, we give the definition of the likelihood ratio order.

**Definition 2.** *Let X and Y be two discrete non negative random variables with $l_X = l_Y = 0$. We say that:*

*(i) X is less than Y in the stochastic order, denoted by $X \leq_{st} Y$, if*

$$S_X(x) \leq S_Y(x), \text{ for all } x \in \mathbb{N}_0.$$

*(ii) X is less than Y in the hazard rate order, denoted by $X \leq_{hr} Y$, if*

$$h_X(x) \geq h_Y(x), \text{ for all } x \in \mathbb{N}_0 \text{ such that } x \leq k^*,$$

*where $k^* = \max\{u_X, u_Y\}$.*

*(iii) X is less than Y in the mean residual life order, denoted by $X \leq_{mrl} Y$, if*

$$m_X(x) \leq m_Y(x), \text{ for all } x \in \mathbb{N}_0 \cup \{-1\} \text{ such that } x \leq \min\{u_X, u_Y\}.$$

*(iv) X is less than Y in the likelihood ratio order, denoted by $X \leq_{lr} Y$, if*

$$\frac{f_X(x)}{f_Y(x)} \text{ is decreasing over the union of the supports,}$$

*where $a/0$ is taken to be equal to $+\infty$ whenever $a > 0$.*

**Remark 2.** *Taking into account that $S_X(x+1) = S_X(x) - f_X(x)$, it is easy to see that the condition in Definition 2 (ii) is equivalent to the following one*

$$\frac{S_X(x)}{S_Y(x)} \text{ is decreasing over } \{0, \ldots, k^*\}. \quad (1)$$

*The ratio $\frac{f_X(x)}{f_Y(x)}$ is usually known as the likelihood ratio.*

It is well-known that the following implications between the previous criteria hold

$$\begin{array}{ccc}
X \leq_{lr} Y \implies & X \leq_{hr} Y & \implies & X \leq_{st} Y \\
& \Downarrow & & \Downarrow \\
& X \leq_{mrl} Y & \implies & E[X] \leq E[Y].
\end{array} \quad (2)$$

## 3. Sufficient Conditions Based on the Unimodality of the Likelihood Ratio

According to (2), the monotonicity of the likelihood ratio is a sufficient condition for the hazard rate order (and the mean residual life and stochastic orders). However, there are situations where the likelihood ratio is non monotone. The purpose of this section is to provide sufficient conditions for the hazard rate and mean residual life orders when the likelihood ratio order is not satisfied, or equivalently, when the likelihood ratio is non monotone. Along the paper, we will denote by $l$ the ratio of the mass probability functions of $X$ and $Y$, that is, the likelihood ratio will be denoted by

$$l(x) = \frac{f_X(x)}{f_Y(x)}, \text{ for all } x \in \{0, \ldots, k^*\}.$$

In addition, we will assume that $l$ is (strictly) unimodal, that is, there exists a value $0 < x_0 < k^*$, such that $l(x)$ is increasing (not constant) over $\{0, \ldots, x_0\}$ and decreasing (not constant) over $\{x_0, \ldots, k^*\}$. Despite the fact that a unimodal function can be increasing or decreasing according to the literature, the distinction among strictly unimodal functions and monotone functions is fundamental in our study. Therefore, when we assume unimodality, we mean strict unimodality.

The first theorem states the behavior of the ratio of the survival functions when $l(x)$ is unimodal. As we will see, under some additional condition, the hazard rate order is satisfied when the likelihood ratio order does not hold. As we have mentioned previously, along this section we will assume that $l_X = l_Y = 0$.

**Theorem 1.** *Let X and Y be two discrete random variables with supports $Supp(X) \subseteq Supp(Y)$. If l is unimodal over $Supp(Y)$, then*

(i)  $S_X(x)/S_Y(x)$ *is decreasing over* $\{0, 1, \ldots, k^*\}$ *if, and only if,* $l(0) \geq 1$.
(ii) $S_X(x)/S_Y(x)$ *is unimodal over* $\{0, 1, \ldots, k^*\}$ *if, and only if,* $l(0) < 1$.

**Proof.** Let us assume that $l(x)$ is increasing for all $0 \leq x \leq k_0$ and decreasing for all $k_0 \leq x \leq k^*$, where $0 < k_0 < k^*$.

Then $[X|X \geq k_0] \leq_{lr} [Y|Y \geq k_0]$ and, from (2), $[X|X \geq k_0] \leq_{hr} [Y|Y \geq k_0]$. According to (1), this inequality is equivalent to the following property

$$\frac{S_X(x)}{S_Y(x)} \text{ is decreasing over } \{k_0, \ldots, k^*\}. \tag{3}$$

Now, let us define the following function

$$t(x) = S_X(x) - l(x)S_Y(x), \text{ for all } x \in \{0, \ldots, k_0\}.$$

Since $l(x)$ is increasing over such set, we have for all $x_1 < x_2 \leq x_0$ that

$$\begin{aligned} t(x_1) &= \sum_{j=x_1}^{x_2-1} (f_X(j) - l(x_1)f_Y(j)) + \sum_{j=x_2}^{\infty} (f_X(j) - l(x_1)f_Y(j)) \\ &\geq \sum_{j=x_2}^{\infty} (f_X(j) - l(x_1)f_Y(j)) \geq \sum_{j=x_2}^{\infty} (f_X(j) - l(x_2)f_Y(j)) \\ &= t(x_2). \end{aligned}$$

Therefore, $t(x)$ is decreasing for all $x \in \{0, \ldots, k_0\}$. Furthermore, since $l(x)$ is decreasing for all $x \geq k_0$, it is easy to see that $t(k_0) \leq 0$.

Now, let us distinguish the following cases.

(i) $l(0) \geq 1$: In this case, we have that

$$t(0) = S_X(0) - l(0)S_Y(0)) \leq (S_X(0) - S_Y(0)) = 0$$

and, therefore, $t(x) \leq 0$ for all $x \in \{0, \ldots, k_0\}$. Therefore,

$$\frac{S_X(x)}{S_Y(x)} \leq l(x) \text{ for all } x \in \{0, \ldots, k_0\},$$

or, equivalently, $h_X(x) \geq h_Y(x)$, for all $x \in \{0, \ldots, k_0\}$, which is equivalent to $S_X(x)/S_Y(x)$ be decreasing over $\{0, \ldots, k_0\}$. Therefore, $S_X(x)/S_Y(x)$ is decreasing over $\{0, \ldots, k^*\}$. It is easy to see that, if $S_X(x)/S_Y(x)$ is decreasing over $\{0, \ldots, k^*\}$, then $l(0) \geq 1$.

(ii) $l(0) < 1$: In this case, we have that

$$t(0) = S_X(0) - l(0)S_Y(0) > S_X(0) - S_Y(0) = 0.$$

Then, there exists a value $0 < k_m \leq k_0$ such that $t(x) > 0$ for all $0 \leq x < k_m$ and $t(x) \leq 0$ for all $k_m \leq x \leq k_0$. Analogously to the previous case, we get that $S_X(x)/S_Y(x)$ is (strictly) increasing for all $0 \leq x \leq k_m$ and decreasing for all $k_m \leq x \leq k_0$. Therefore, $S_X(x)/S_Y(x)$ is unimodal over $\{0, 1, \ldots, k_0\}$. Given the characterization in (i), if $S_X(x)/S_Y(x)$ is unimodal over $\{0, 1, \ldots, k^*\}$, then $l(0) < 1$.

To sum up, combining (3) with the previous discussion, we conclude that $S_X/S_Y$ is decreasing if, and only if, $l(0) \geq 1$ and it is unimodal otherwise. □

We observe that the conditions on $l(0)$ can be described in terms of the comparison of probability of failed components before to put into operation. For example if $X$ is the lifetime of a device after a burn-in process then $l(0) = 0 < 1$.

Therefore, there are situations where $l(x)$ is unimodal but the hazard rate order holds and, consequently, the mean residual life order too. Now, a natural question arises: In the situations where the ratio of the survival functions is unimodal, does the mean residual life order hold? Ref. [10] provides a result in this direction in the continuous case. Given that the mean residual life function in the discrete case can be obtained from the formula in the continuous case, the next result trivially follows from Theorem 2.3 in [10].

**Theorem 2.** *Let $X$ and $Y$ be two discrete random variables with supports $Supp(X) \subseteq Supp(Y)$, and survival functions such that $S_X(x)/S_Y(x)$ is unimodal over $\{0, \ldots, k^*\}$. Then, $X \leq_{mrl} Y$ if, and only if, $E[X] \leq E[Y]$.*

**Remark 3.** *We also want to point out that the mean residual life order holds when the survival ratio $S_X/S_Y$ is unimodal and $E[X] \leq E[Y]$, but the stochastic order does not hold. Let us justify this assertion. Given that $S_X/S_Y$ is unimodal we have that $S_X(x) \geq S_Y(x)$, for all $0 \leq x \leq k_0$. Furthermore, there exists at least a value $x'$ such that $S_X(x') < S_Y(x')$, since otherwise we have $E[X] > E[Y]$. Therefore, the stochastic order does not hold in such situations in any sense.*

Taking into account this remark and Theorems 1 and 2, we can state the following corollary that characterizes the hazard rate and the mean residual life orders when the ratio of the mass probability functions is unimodal.

**Corollary 1.** *Let $X$ and $Y$ be two discrete random variables with supports $Supp(X) \subseteq Supp(Y) \subseteq \mathbb{N}_0$. If $l$ is unimodal over $Supp(Y)$, then*

(i) *$X \leq_{hr} Y$ if, and only if, $l(0) \geq 1$.*
(ii) *$X \leq_{mrl} Y$ ($X \nleq_{st} Y$) if, and only if, $l(0) < 1$ and $E[X] \leq E[Y]$.*

As we have mentioned in the introduction, ref. [11] has dealt with the unimodality of $l$ to establish conditions for the usual stochastic order for two discrete random variables. In particular, they prove that the umimodality and the condition in Corollary 1 (i) imply the stochastic order. Therefore, we have improved their result, concluding that the hazard rate order holds under such condition, which is stronger than the stochastic one. In fact, they also ask for another assumption on the limit of $l$ which is not necessary. Furthermore, we also state what happens in the situations where the limit condition is not satisfied, which characterizes all the possible scenarios in terms of these stochastic orders when the likelihood ratio *l* is unimodal.

## 4. Applications

In this section we apply Corollary 1 and Theorem 2 to compare two random variables belonging to some of the most important discrete parametric families in reliability. The starting point is to apply Corollary 1 to look for situations where the likelihood ratio is unimodal. Unluckily, there are situations where this is not an easy task. It is well-known that a sufficient condition for the unimodality is the property of logconcavity. Let us formally define it.

**Definition 3.** *Let $A = \{a_j\}_{j=0}^{k}$ be a succession of real values, then A is said to be logconcave if $a_j^2 \geq a_{j+1}a_{j-1}$, for all $1 \leq j < k$, where k can be finite or infinite.*

It is easy to see that if $A = \{a_j\}_{j=0}^{k}$ is a succession of positive real values, then

$$\frac{a_n}{a_{n-1}} \geq \frac{a_{n+1}}{a_n}, \text{ for all } n < k. \quad (4)$$

A well-known property is that if $A = \{a_j\}_{j=0}^{k}$ is logconcave and positive then *A* is (not strictly) unimodal, that is, there exists a value $j \in \mathbb{N}_0$ such that $a_0 \leq \ldots \leq a_{j-1} \leq a_j \geq a_{j+1} \geq \ldots \geq a_n \geq \ldots$ (see [13]) and, therefore, *A* can be increasing, decreasing, or unimodal. In order to apply the results of the previous section, it is fundamental to determine if the succession of the likelihood ratio values is monotone or strictly unimodal. Next, we give a result that characterizes these two situations.

**Proposition 1.** *Let $A = \{a_j\}_{j=0}^{k}$ be a logconcave and positive succession, and let us assume that $\lim_{j \to k} a_{j+1}/a_j$ exists. Then,*

(i)　*A is decreasing if, and only if, $a_1/a_0 \leq 1$.*
(ii)　*A is increasing if, and only if, $\lim_{j \to k} a_{j+1}/a_j \geq 1$.*
(iii)　*A is unimodal if, and only if, $a_1/a_0 > 1$ and $\lim_{j \to k} a_{j+1}/a_j < 1$.*

**Proof.** Case (i): Let us assume first that *A* is decreasing, then we have that $a_1/a_0 \leq 1$. Let us prove now the reversed implication. From (4) we have that $a_{j+1}/a_j$ is decreasing, therefore $a_{j+1}/a_j \leq a_1/a_0 \leq 1$ and then *A* is decreasing.

Case (ii): It follows under similar arguments as in case (i).

Case (iii): If *A* is unimodal the result follows from (i) and (ii). □

Next, we combine the previous proposition with Corollary 1 to provide sufficient conditions for the likelihood, hazard rate, and mean residual life orders, under the assumption of logconcavity of the likelihood ratio.

**Proposition 2.** *Let X and Y be two discrete random variables with support in $\mathbb{N}_0$, such that $Supp(X) \subseteq Supp(Y)$, and mass probability functions such that $l(x)$ is logconcave over $0 \leq x \leq u_X$. Then we have the following results:*

(i)　$X \leq_{lr} Y$ *if, and only if, $l(1)/l(0) \leq 1$.*
(ii)　*If*

$$u_X < u_Y \text{ and } \frac{l(1)}{l(0)} > 1,$$

*or*

$$u_X = u_Y, \frac{l(1)}{l(0)} > 1 \text{ and } \lim_{j \to u_X} \frac{l(j+1)}{l(j)} < 1,$$

*then*

(ii.a) $X \leq_{hr} Y$ ($X \not\leq_{lr} Y$) if, and only if, $l(0) \geq 1$.
(ii.b) $X \leq_{mrl} Y$ ($X \not\leq_{hr} Y$) if, and only if, $l(0) < 1$ and $E[X] \leq E[Y]$.

Now, we provide an application of this result to the most classical parametric discrete models of interest in reliability, that is, the binomial, negative binomial and Poisson distributions, which are particular cases of the Panjer's family. It is easy to see that $l(x)$ is logconcave for any pair of Panjer's distributions and, therefore, we can apply Proposition 2 to compare two random variables belonging to each one of these families. Let us formally define the Panjer's family.

**Example 1.** *Panjer's family:* Ref. [14] *considered a family of discrete random variables with mass probability function satisfying the following relationship*

$$f_X(x) = \frac{ax+b}{x} f_X(x-1), \text{ for all } x \in \mathbb{N},$$

*where $a, b \in \mathbb{R}$ and $f(0) > 0$. Ref. [15] proved that the only three discrete distributions that satisfy the previous relationship are the Poisson, binomial, and negative binomial distributions, with the following parameters:*

(a) *Poisson $P(\lambda)$: $a = 0$, $b = \lambda$ and $f_X(0) = e^{-\lambda}$, where $\lambda \in \mathbb{R}_+$.*
(b) *Binomial $B(n, p)$: $a = -p/(1-p)$, $b = (n+1)p/(1-p)$ and $f_X(0) = (1-p)^n$, where $n \in \mathbb{N}$ and $p \in (0,1)$.*
(c) *Negative binomial $NB(p, r)$: $a = p - 1$, $b = (r-1)(p-1)$ and $f_X(0) = p^r$, where $r \in \mathbb{N}$ and $p \in (0,1)$.*

*Let $X$ and $Y$ be defined according to the Panjer's family with parameters $a_1$, $b_1$ and $a_2$, $b_2$, respectively, $f_X(0), f_Y(0) > 0$ and $u_X \leq u_Y$. Since $l(x) = \frac{a_1 x + b_1}{a_2 x + b_2} l(x-1) > 0$, for all $x \in \{0, 1, \ldots, u_X\}$, we have that*

$$\frac{l(x+1)}{l(x)} = \frac{a_1 x + b_1}{a_2 x + b_2}, \text{ for all } x \in \{0, 1, \ldots, u_X\},$$

*and it is easy to see that*

$$\frac{l(x+1)}{l(x)} \text{ is decreasing if, and only if, } a_1 b_2 \leq a_2 b_1.$$

*Therefore, according to Proposition 2, the succession $l(x)$ is logconcave for $x$ in $\{0, \ldots, u_X\}$ if, and only if, $a_1 b_2 \leq a_2 b_1$. In addition, we have that*

$$\frac{l(1)}{l(0)} = \frac{a_1 + b_1}{a_2 + b_2}.$$

*Then, by the previous discussion, we can easily state the following results applying Proposition 2.*

(a) **Poisson:** *Let $X \sim P(\lambda_1)$ and $Y \sim P(\lambda_2)$. Then, $X \leq_{lr} Y$ if, and only if,*

$$\lambda_1 \leq \lambda_2.$$

(b) **Binomial:** *Let $X \sim B(n_1, p_1)$ and $Y \sim B(n_2, p_2)$, where $n_1 \leq n_2$ without loss of generality. Then,*

(b.1) *$X \leq_{lr} Y$ if, and only if,*

$$n_1 p_1 (1 - p_2) \leq n_2 p_2 (1 - p_1).$$

(b.2) *$X \leq_{hr} Y$ (and $X \not\leq_{lr} Y$) if, and only if,*

$$\left\{ \begin{array}{rcl} n_1 p_1 (1 - p_2) & > & n_2 p_2 (1 - p_1) \\ (1 - p_1)^{n_1} & \geq & (1 - p_2)^{n_2}. \end{array} \right\}$$

(b.3) $X \leq_{mrl} Y$ (and $X \not\leq_{hr} Y$) if, and only if,

$$\left\{ \begin{array}{rcl} n_1 p_1(1-p_2) & > & n_2 p_2(1-p_1) \\ (1-p_1)^{n_1} & < & (1-p_2)^{n_2} \\ n_1 p_1 & \leq & n_2 p_2. \end{array} \right\}$$

(c) **Negative binomial:** Let $X \sim NB(p_1, r_1)$ and $Y \sim NB(p_2, r_2)$, where $r_1 \geq r_2$ without loss of generality. Then,

(c.1) $X \leq_{lr} Y$ if, and only if,

$$r_1(1-p_1) \leq r_2(1-p_2).$$

(c.2) $X \leq_{hr} Y$ (and $X \not\leq_{lr} Y$) if, and only if,

$$\left\{ \begin{array}{rcl} r_1(1-p_1) & > & r_2(1-p_2) \\ p_1 & > & p_2 \\ p_1^{r_1} & \geq & p_2^{r_2}. \end{array} \right\}$$

(c.3) $X \leq_{mrl} Y$ (and $X \not\leq_{hr} Y$) if, and only if,

$$\left\{ \begin{array}{rcl} r_1(1-p_1) & > & r_2(1-p_2) \\ p_1 & > & p_2 \\ p_1^{r_1} & < & p_2^{r_2} \\ r_1(1-p_1)p_2 & \leq & r_2(1-p_2)p_1. \end{array} \right\}$$

Let us now make some comments on the previous results. As for the Poisson family regards, it is an already well-known result in the literature. As the binomial is concerned, the previous result improves Theorem 1.1 in [11] since they establish the stochastic order under the same conditions as we do for the hazard rate order. Finally, regarding the binomial negative family, our result also improves Theorem 1.1 in [11] because in those cases where $r_1 \geq r_2$ the likelihood order or the hazard rate order hold under the same conditions as the stochastic order. Let us point out that the assumption $r_1 \geq r_2$ is not a restriction since, otherwise, the roles of $X$ and $Y$ can be exchanged.

Next, we consider a discrete version of the continuous Weibull distribution. Ref. [16] introduced the Weibull discrete family in the context of reliability in discrete time. The problem of estimation of its parameters can be seen in [17,18]. In this case the results follow from the analysis of the ratio of the survival functions.

**Example 2.** *Discrete Weibull distribution:* Let $X$ be a discrete random variable. We say that $X$ follows a discrete **Weibull** distribution with parameters $\alpha > 0$ and $\beta > 0$, denoted by $X \sim W(\alpha, \beta)$, if its mass probability function is given by

$$f_X(x) = q^{x^\beta} - q^{(x+1)^\beta}, \text{ for all } x \in \mathbb{N}_0,$$

where $q = \exp\{-1/\alpha^\beta\} \in (0,1)$. The survival function has also a closed expression and it is given by

$$S_X(x) = q^{x^\beta}, \text{ for all } x \in \mathbb{N}_0.$$

Let $X \sim W(\alpha_1, \beta_1)$ and $Y \sim W(\alpha_2, \beta_2)$ be two random variables with Weibull distributions and let us assume without loss of generality that $\beta_1 \geq \beta_2$. Then,

(a) $X \leq_{hr} Y$ if, and only if.

$$\alpha_1^{\beta_1} \leq \alpha_2^{\beta_2}. \tag{5}$$

(b) $X \leq_{mrl} Y$ (and $X \not\leq_{hr} Y$) if, and only if,

$$\left\{ \begin{array}{rcl} \alpha_1^{\beta_1} & > & \alpha_2^{\beta_2} \\ E[X] & \leq & E[Y]. \end{array} \right\} \tag{6}$$

*Let us see the justification of the previous results. Denoting by $L(x)$ the ratio of the survival functions, that is, $L(x) = S_X(x)/S_Y(x) = q_1^{x^{\beta_1}}/q_2^{x^{\beta_2}}$, it is not difficult to prove that $L$ is unimodal over $\mathbb{R}^+$ (not necessarily over $\mathbb{N}_0$), where the maximum is attained at*

$$x_0 = \left(\beta_1 \alpha_2^{\beta_2}/(\beta_2 \alpha_1^{\beta_1})\right)^{1/(\beta_2 - \beta_1)}.$$

*Given that $L(0) = 1$, it is easy to see that $L(x)$ is decreasing for all $x \in \mathbb{N}_0$ if, and only if, $L(1) \leq 1$ and it is unimodal otherwise. In other words, we have that $X \leq_{hr} Y$ if, and only if, (5) holds. Furthermore, if (5) does not hold, we have that $L(x)$ is unimodal and, applying Theorem 2, we have that $X \leq_{mrl} Y$ if, and only if, (6) holds.*

*In Figure 1, we provide examples of the two previous situations.*

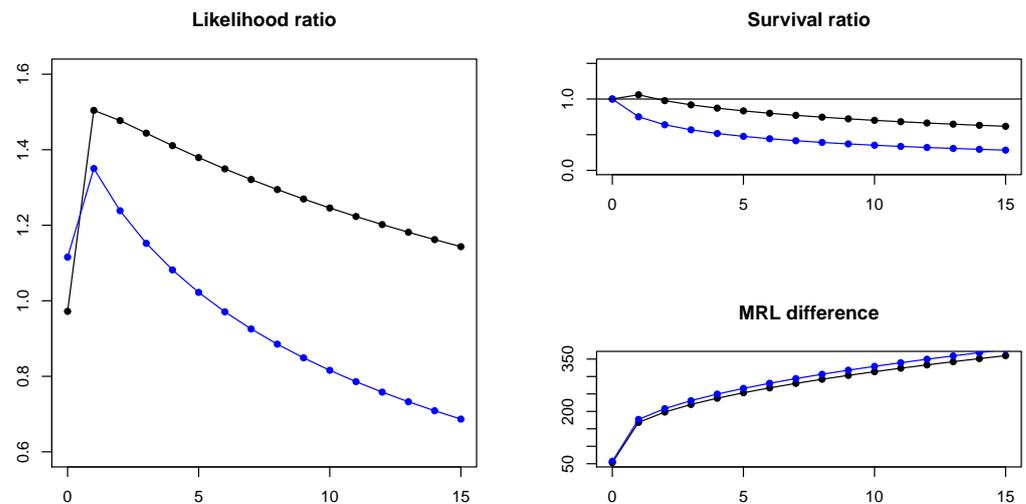

**Figure 1.** Likelihood ratio, ratio of survival functions and mean residual life difference for case (a) where $X \sim W(0.3, 0.3)$ and $Y \sim W(0.5, 0.2)$ (black) and case (b) where $X \sim W(0.75, 0.3)$ and $Y \sim W(0.5, 0.2)$ (blue).

The upcoming applications are based on the unimodality of the likelihood ratio. First, we consider a generalization of the Poisson family that covers some other situations apart from the Poisson distribution. A comprehensive study of this distribution is available in [19].

**Example 3.** *Generalized Poisson' family:* Let $X$ be a discrete random variable. We say that $X$ *follows a **generalized Poisson** distribution with parameters $\theta > 0$ and $0 < \lambda < 1$, denoted by $X \sim P(\theta, \lambda)$, if its mass probability function is given by*

$$f_X(x) = \frac{\theta(\theta + x\lambda)^{x-1}}{x!}\exp\{-(\theta + \lambda x)\}, \text{ for all } x \in \mathbb{N}_0.$$

*The distribution function has no closed expression and the mean is given by $\theta/(1-\lambda)$.*

*Let $X \sim GP(\theta_1, \lambda_1)$ and $Y \sim GP(\theta_2, \lambda_2)$ be two random variables with generalized Poisson distributions where $\lambda_1 \leq \lambda_2$ without loss of generality. Then,*

(a)　$X \leq_{lr} Y$ *if, and only if,*

$$\lambda_2 \theta_1 < \lambda_1 \theta_2 \tag{7}$$

*or*

$$\left\{\begin{array}{rcl} \lambda_2 \theta_1 & > & \lambda_1 \theta_2 \\ \log\left(\frac{\theta_2}{\theta_1}\right) & \geq & \frac{\lambda_2 \theta_1 - \lambda_1 \theta_2}{\theta_1 \theta_2} + \lambda_2 - \lambda_1. \end{array}\right\} \tag{8}$$

(b)　$X \leq_{hr} Y$ ( and $X \not\leq_{lr} Y$) if, and only if,

$$\left\{ \begin{array}{rcl} \lambda_2 \theta_1 & > & \lambda_1 \theta_2 \\ \log\left(\frac{\theta_2}{\theta_1}\right) & < & \frac{\lambda_2 \theta_1 - \lambda_1 \theta_2}{\theta_1 \theta_2} + \lambda_2 - \lambda_1 \\ \theta_1 & \leq & \theta_2. \end{array} \right\} \quad (9)$$

(c)　$X \leq_{mrl} Y$ ( and $X \not\leq_{hr} Y$) if, and only if,

$$\left\{ \begin{array}{rcl} \lambda_2 \theta_1 & > & \lambda_1 \theta_2 \\ \log\left(\frac{\theta_2}{\theta_1}\right) & < & \frac{\lambda_2 \theta_1 - \lambda_1 \theta_2}{\theta_1 \theta_2} + \lambda_2 - \lambda_1 \\ \theta_1 & > & \theta_2 \\ \theta_1(1 - \lambda_2) & \leq & \theta_2(1 - \lambda_1). \end{array} \right\} \quad (10)$$

*Let us justify the previous results. First, we show that*

$$l(x) = \frac{\theta_1}{\theta_2}\left(\frac{\theta_1 + x\lambda_1}{\theta_2 + x\lambda_2}\right)^{x-1} \exp\{-\theta_1 + \theta_2 + x(\lambda_1 - \lambda_2)\}$$

*is decreasing or unimodal over* $\mathbb{N}_0$ *whenever* $\lambda_1 \leq \lambda_2$. *In order to prove it, we consider the behavior of* $l(x)$ *over* $\mathbb{R}_+$ *or, equivalently, the behavior of* $\log l(x)$ *over* $\mathbb{R}_+$. *After some computations, we get that* $(\log l(x))' = c_1(x) - c_2(x)$ *where*

$$c_1(x) = \left[\frac{(x-1)(\theta_2\lambda_1 - \theta_1\lambda_2)}{(\theta_1 + x\lambda_1)(\theta_2 + x\lambda_2)} + \lambda_2 - \lambda_1\right]$$

*and*

$$c_2(x) = \left[\log\left(\frac{\theta_2 + x\lambda_2}{\theta_1 + x\lambda_1}\right)\right].$$

*On the one hand, it is easy to see that* $c_2(x)$ *is increasing (decreasing), whenever* $\lambda_2\theta_1 \geq (\leq)\lambda_1\theta_2$. *We also have that*

$$\lim_{x \to 0^+} c_2(x) = \log\left(\frac{\theta_2}{\theta_1}\right) \text{ and } \lim_{x \to +\infty} c_2(x) = \log\left(\frac{\lambda_2}{\lambda_1}\right).$$

*On the other hand, we get that* $c_1(x)$ *is decreasing for all* $0 \leq x \leq x_0$ *and increasing for all* $x \geq x_0$, *where*

$$x_0 = 1 + \sqrt{1 + \frac{\theta_2/\theta_1 + \lambda_1\theta_2 + \lambda_2\theta_1}{\lambda_2/\lambda_1}}.$$

*We also have that*

$$\lim_{x \to 0^+} c_1(x) = \frac{\lambda_2\theta_1 - \lambda_1\theta_2}{\theta_1\theta_2} + \lambda_2 - \lambda_1 \text{ and } \lim_{x \to +\infty} c_1(x) = \lambda_2 - \lambda_1.$$

*Taking into account the previous information, we get that the modes of* $l(x)$ *are the crossing points of an increasing (decreasing) function with a parabola with a minimum. Now, since* $\lambda_1 \leq \lambda_2$ *implies*

$$\lim_{x \to +\infty} c_1(x) = \lambda_2 - \lambda_1 \leq \log\left(\frac{\lambda_2}{\lambda_1}\right) = \lim_{x \to +\infty} c_2(x),$$

*it is easy to see that* $l(x)$ *is decreasing for all* $x \in \mathbb{N}_0$ *if, and only if,* (7) *or* (8) *holds. Otherwise,* $l(x)$ *is unimodal. Now, applying Corollary 1 and taking into account that*

$$\lim_{x \to 0^+} l(x) = \exp\{\theta_2 - \theta_1\},$$

*we get that* $X \leq_{hr} Y$ *(and* $X \not\leq_{lr} Y$*) if, and only if,* (9) *holds and* $X \leq_{mrl} Y$ *(and* $X \not\leq_{hr} Y$*) if, and only if,* (10) *holds.*

*In Figure 2, we provide examples of the previous situations.*

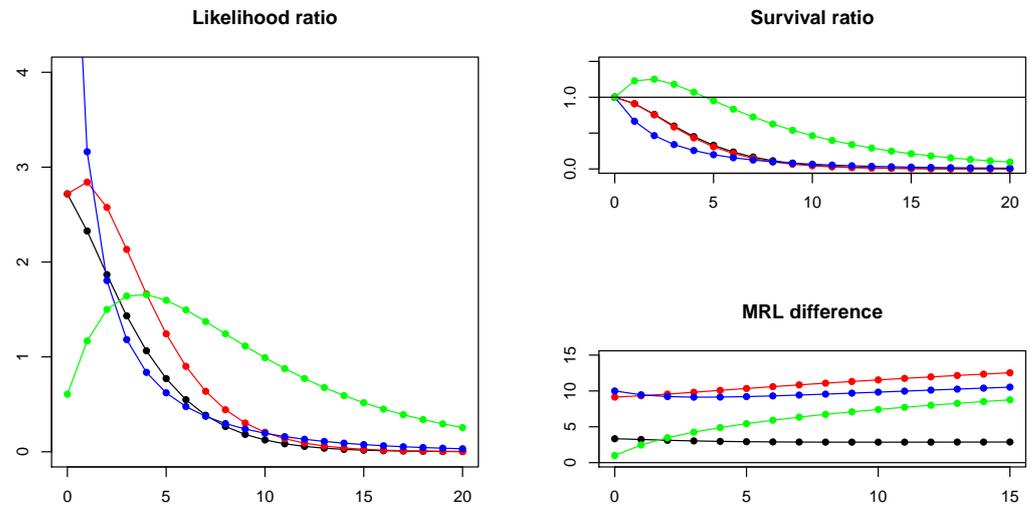

**Figure 2.** Likelihood ratio, ratio of survival functions, and mean residual life difference for case (a) (condition (7)) where $X \sim GP(2, 0.25)$ and $Y \sim GP(3, 0.5)$ (black), case (a) (conditions (8)) where $X \sim GP(1, 0.5)$ and $Y \sim GP(3, 0.75)$ (blue), case (b) where $X \sim GP(2, 0.3)$ and $Y \sim GP(3, 0.75)$ (red) and case (c) where $X \sim GP(2, .0.5)$ and $Y \sim GP(1, 0.75)$ (green).

Finally, we deal with the Hurwitz-Lerch family. Ref. [20] introduced this family as a unified representation of some discrete distributions. Particular cases of this family are Zipf, Zipf-Mandelbrot, and Good distributions. Ref. [21] study several properties of this family as well as their application in reliability. The reader can find a detailed study and further references in [12]. As we will see, the left extreme of the support is 1 instead of 0. The results provided in this paper have been given for random variables with left extreme equal to 0, but it is clear that can be used for non zero common left extremes considering a change of origin.

**Example 4.** *Hurwitz-Lerch's family:* Let X be a discrete random variable. We say that X follows a **Hurwitz-Lerch** distribution with parameters $s \geq 0$, $0 < z \leq 1$ and $0 \leq a \leq 1$, denoted by $X \sim HL(z, s, a)$, if its mass probability function is given by

$$f_X(x) = \frac{z^x}{T(z, s, a)(a + x)^{s+1}}, \text{ for all } x \in \mathbb{N},$$

where $T(z, s, a) = \sum_{x=1}^{\infty} z^x/(a+x)^{s+1}$. Once again, the distribution function does not have a closed expression and the mean is given by $T(z, s-1, a)/T(z, s, a) - a$.

Let $X \sim HL(z_1, s_1, a_1)$ and $Y \sim HL(z_2, s_2, a_2)$ be two random variables with Hurwitz-Lerch distributions such that $a_1 \geq a_2$ without loss of generality. Then,

(a)  If $a_1 = a_2 = a$, then

(a.1) $X \leq_{lr} Y$ if, and only if,

$$\left\{ \begin{array}{rcl} z_1 & \leq & z_2 \\ l(1) & \geq & l(2). \end{array} \right\}$$

(a.2) $X \leq_{hr} Y$ (and $X \not\leq_{lr} Y$) if, and only if,

$$\left\{ \begin{array}{rcl} z_1 & \leq & z_2 \\ l(1) & < & l(2) \\ \frac{T(z_2, s_2, a)(a+1)^{s_2 - s_1}}{z_2} & \geq & \frac{T(z_1, s_1, a)}{z_1}. \end{array} \right\}$$

(a.3) $X \leq_{mrl} Y$ (and $X \nleq_{hr} Y$) if, and only if

$$\left\{ \begin{array}{rcl} z_1 & \leq & z_2 \\ l(1) & < & l(2) \\ \frac{T(z_2,s_2,a)(a+1)^{s_2-s_1}}{z_2} & < & \frac{T(z_1,s_1,a)}{z_1} \\ E[X] & \leq & E[Y]. \end{array} \right\}$$

(b) If $a_1 > a_2$, then

(b.1) $X \leq_{lr} Y$ if, and only if, either (11), or (12) or (13) holds, where

$$\left\{ \begin{array}{rcl} z_1 & \leq & z_2 \\ \frac{s_1+1}{a_1+1} + \log\left(\frac{z_2}{z_1}\right) & \geq & \frac{s_2+1}{a_2+1}. \end{array} \right\} \quad (11)$$

$$\left\{ \begin{array}{rcl} z_1 & < & z_2 \\ \frac{s_1+1}{a_1+1} + \log\left(\frac{z_2}{z_1}\right) & < & \frac{s_2+1}{a_2+1} \\ l(1) & \geq & l(2). \end{array} \right\} \quad (12)$$

$$\left\{ \begin{array}{rcl} z_1 & = & z_2 \\ s_1 & > & s_2 \\ \frac{s_1+1}{a_1+1} & < & \frac{s_2+1}{a_2+1} \\ l(1) & \geq & l(2). \end{array} \right\} \quad (13)$$

(b.2) $X \leq_{hr} Y$ (and $X \nleq_{lr} Y$) if, and only if, either (14) or (15) holds, where

$$\left\{ \begin{array}{rcl} z_1 & < & z_2 \\ \frac{s_1+1}{a_1+1} + \log\left(\frac{z_2}{z_1}\right) & < & \frac{s_2+1}{a_2+1} \\ l(1) & < & l(2) \\ \frac{T(z_2,s_2,a_2)(a_2+1)^{s_2+1}}{z_2} & \geq & \frac{T(z_1,s_1,a_1)(a_1+1)^{s_1+1}}{z_1}. \end{array} \right\} \quad (14)$$

$$\left\{ \begin{array}{rcl} z_1 & = & z_2 \\ s_1 & > & s_2 \\ (s_1+1)(a_2+1) & < & (s_2+1)(a_1+1) \\ l(1) & < & l(2) \\ \frac{T(z_2,s_2,a_2)(a_2+1)^{s_2+1}}{z_2} & \geq & \frac{T(z_1,s_1,a_1)(a_1+1)^{s_1+1}}{z_1}. \end{array} \right\} \quad (15)$$

(b.3) $X \leq_{mrl} Y$ (and $X \nleq_{hr} Y$) if, and only if, either (16) or (17) holds, where

$$\left\{ \begin{array}{rcl} z_1 & < & z_2 \\ \frac{s_1+1}{a_1+1} + \log\left(\frac{z_2}{z_1}\right) & < & \frac{s_2+1}{a_2+1} \\ l(1) & < & l(2) \\ \frac{T(z_2,s_2,a_2)(a_2+1)^{s_2+1}}{z_2} & < & \frac{T(z_1,s_1,a_1)(a_1+1)^{s_1+1}}{z_1} \\ E[X] & \leq & E[Y]. \end{array} \right\} \quad (16)$$

$$\left\{ \begin{array}{rcl} z_1 & = & z_2 \\ s_1 & > & s_2 \\ (s_1+1)(a_2+1) & < & (s_2+1)(a_1+1) \\ l(1) & < & l(2) \\ \frac{T(z_2,s_2,a_2)(a_2+1)^{s_2+1}}{z_2} & < & \frac{T(z_1,s_1,a_1)(a_1+1)^{s_1+1}}{z_1} \\ E[X] & \leq & E[Y]. \end{array} \right\} \quad (17)$$

In Figure 3, we provide some examples of the previous situations.

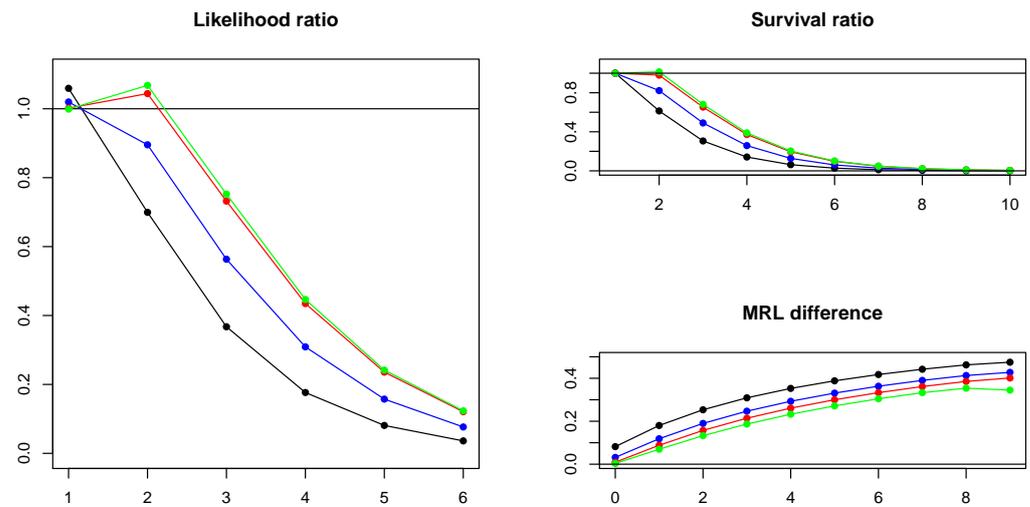

**Figure 3.** Likelihood ratio, ratio of survival functions, and mean residual life difference for case (a.1) (conditions (11)) where $X \sim HL(0.2, 1, 0.7)$ and $Y \sim HL(0.5, 1.5, 0.3)$ (black), case (a.1) (conditions (12)) where $X \sim HL(0.2, 1, 0.7)$ and $Y \sim HL(0.5, 2, 0.3)$ (blue), case (b.2) (condition (14)) where $X \sim HL(0.2, 1, 0.7)$ and $Y \sim HL(0.5, 2.3, 0.3)$ (red) and case (b.3) (condition (16)) where $X \sim HL(0.2, 1.5, 0.7)$ and $Y \sim HL(0.5, 2.75, 0.3)$ (green).

*Additionally, we can provide similar results on the different stochastic orders in the reverse sense that do not follow just exchanging the role of the random variables.*

(a) $X \geq_{lr} Y$ if, and only if, one of the following sets of conditions holds:

$$\left\{\begin{array}{rcl} z_1 & = & z_2 \\ s_1 & \leq & s_2. \end{array}\right\}$$

$$\left\{\begin{array}{rcl} z_1 & > & z_2 \\ \frac{s_1+1}{a_1+1} + \log\left(\frac{z_2}{z_1}\right) & \leq & \frac{s_2+1}{a_2+1}. \end{array}\right\}$$

$$\left\{\begin{array}{rcl} z_1 & > & z_2 \\ \frac{s_1+1}{a_1+1} + \log\left(\frac{z_2}{z_1}\right) & > & \frac{s_2+1}{a_2+1} \\ l(1) & \leq & l(2). \end{array}\right\}$$

(b) $X \geq_{hr} Y$ (and $X \not\geq_{lr} Y$) if, and only if,

$$\left\{\begin{array}{rcl} z_1 & > & z_2 \\ \frac{s_1+1}{a_1+1} + \log\left(\frac{z_2}{z_1}\right) & > & \frac{s_2+1}{a_2+1} \\ l(1) & > & l(2) \\ \frac{T(z_2,s_2,a_2)(a_2+1)^{s_2+1}}{z_2} & \leq & \frac{T(z_1,s_1,a_1)(a_1+1)^{s_1+1}}{z_1}. \end{array}\right\}$$

(c) $X \geq_{mrl} Y$ (and $X \not\geq_{hr} Y$) if, and only if,

$$\left\{\begin{array}{rcl} z_1 & > & z_2 \\ \frac{s_1+1}{a_1+1} + \log\left(\frac{z_2}{z_1}\right) & > & \frac{s_2+1}{a_2+1} \\ l(1) & > & l(2) \\ \frac{T(z_2,s_2,a_2)(a_2+1)^{s_2+1}}{z_2} & > & \frac{T(z_1,s_1,a_1)(a_1+1)^{s_1+1}}{z_1} \\ E[X] & \geq & E[Y]. \end{array}\right\}$$

*Let us provide a justification of these results. In this case, we have that*

$$l(x) = \frac{T(z_2, s_2, a_2)}{T(z_1, s_1, a_1)} \left(\frac{z_1}{z_2}\right)^x \frac{(a_2 + x)^{s_2+1}}{(a_1 + x)^{s_1+1}}.$$

*The study of the changes of monotonicity of $l(x)$ over $\mathbb{R}^+$ is equivalent to study the sign changes of*

$$(s_2 + 1)\frac{a_1 + x}{a_2 + x} - \left(\log\left(\frac{z_2}{z_1}\right)x + a_1 \log\left(\frac{z_2}{z_1}\right) + s_1 + 1\right).$$

*Taking into account that the previous expression is the difference between a convex decreasing function (since $a_1 \geq a_2$) and a linear function, using geometric arguments we conclude that $l(x)$ is either unimodal or monotone, setting the situations where each one occurs. Then, the arguments for the comparison of two members of this family are similar to those provided in the previous example and will not be reproduced with detail due to their length.*

Let us point out that [22] provides several results for the comparison of two Hurwitz-Lerch distributions in the likelihood ratio order. The conditions that we have stated improved their result for the likelihood ratio order, since we have characterized the situations where this criterion holds in some sense. In addition, we give situations where the hazard rate and the mean residual life orders hold.

## 5. Conclusions

In the present paper we have provided several results on the comparison of two discrete distributions according to the likelihood ratio, the hazard rate and the mean residual life orders. These results fill the gap in the literature since there are only results for comparing discrete distributions in the usual stochastic order. The results have been applied to compare some well-known families of discrete distributions in the context of reliability in discrete time.

**Author Contributions:** Conceptualization, F.B. and C.M.-R.; methodology, F.B., C.M.-R. and M.P.; software, C.M.-R.; validation, F.B. and C.M.-R.; formal analysis, F.B., C.M.-R. and M.P.; investigation, F.B., C.M.-R. and M.P.; resources, F.B., C.M.-R. and M.P.; writing—original draft preparation, F.B. and C.M.-R.; writing—review and editing, F.B. and C.M.-R.; visualization, F.B. and C.M.-R.; supervision, F.B. and C.M.-R.; project administration, F.B.; funding acquisition, F.B. All authors have read and agreed to the published version of the manuscript.

**Funding:** This research was funded by Ministerio de Ciencia e Innovación of Spain under grant PID2019-103971GB-I00/AEI/10.13039/501100011033.

**Institutional Review Board Statement:** Not applicable.

**Informed Consent Statement:** Not applicable.

**Data Availability Statement:** Not applicable.

**Conflicts of Interest:** The authors declare no conflict of interest.